\documentclass[a4paper,12pt,leqno]{amsart}
\usepackage[bulgarian,english]{babel}
\usepackage[cp1251]{inputenc}
\usepackage[T2A]{fontenc}
\usepackage{amsmath}
\usepackage{amssymb}

\newcommand{\ds}{\displaystyle}

\newtheorem{theorem}{Theorem}

\newtheorem{lemma}{Lemma}

\title[On the distribution of $\alpha p$ modulo one]
{On the distribution of $\alpha p$ modulo one for primes $p$ of a
special form}

\author{T.L.Todorova \and D.I.Tolev}

\date{}

\begin{document}
\maketitle

\begin{abstract}
A classical problem in analytic number theory is to study the
distribution of $\alpha p$ modulo 1, where $\alpha$ is irrational
and $p$ runs over the set of primes. We consider the subsequence
generated by the primes $p$ such that $p+2$ is an almost-prime
(the existence of infinitely many such $p$ is another topical
result in prime number theory) and prove that its distribution has
a similar property.


\end{abstract}

\section{Introduction and statements of the result}
The famous prime twins conjecture states that there exist
infinitely many primes $p$ such that $p+2$ is a prime too. This
hypothesis is still unproved but there are many approximation to
it established. One of the most interesting of them is due to Chen
\cite{Chen}. In 1973 he proved that there are infinitely many
primes $p$ for which $p+2=\mathrm{P}_2$. (As usual $\mathrm{P}_r$
denotes an integer with no more than $r$ prime factors, counted
according to multiplicity).

Suppose that we have a problem including primes and let $r\ge 2$
be an integer. Having in mind Chen's result we may consider this
problem with primes $p$, such that $p+2=\mathrm{P}_r.$
We will give several examples.

In 1937, Vinogradov \cite{Vin1} proved that every
sufficiently large odd $n$ can be represented in the form
\begin{equation}\label{reprn}
    p_1+p_2+p_3=n,
\end{equation}
where $p_1,\,p_2,\,p_3$ are primes. In 2000 Peneva \cite{Pen} and Tolev
\cite{Tol1} considered~\eqref{reprn} with primes of the form
specified above. It was established in \cite{Tol1} that if $n$ is
sufficiently large and $n\equiv 3\pmod 6$ then~\eqref{reprn} has a
solution in primes $p_1,\,p_2,\,p_3$ such that
\[
p_1+2=\mathrm{P}_2, \quad p_2+2=\mathrm{P}_5,  \quad
p_3+2=\mathrm{P}_7.
\]

Further, in 1938 Hua \cite{Hua} proved that every sufficiently
large $n \equiv 5 \pmod {24}$ can be represented as
\begin{equation}\label{reprn2}
    p^2_1 + p^2_2 + p^2_3 + p^2_4 + p^2_5 = n,
\end{equation}
where $p_1, \ldots , p_5$ are primes. In 2000 Tolev \cite{Tol3}
proved that every sufficiently large $n \equiv 5 \pmod {24}$ can
be represented in the form~\eqref{reprn2} with primes $p_1, \ldots
, p_5$ such that
\begin{equation*}
    p_1 + 2 = P_2, \,p_2 + 2 = P'_2 ,\, p_3 + 2 = P_5,\, p_4 + 2 =
P'_5 , \,p_5 + 2 = P_7.
\end{equation*}

Finally, in 2004 Green and Tao \cite{GrT1} proved their celebrated
theorem stating that for every natural $k\ge 3$ there are
infinitely many arithmetical progression of $k$ different primes.
Later they established (see~\cite{GrT2}) that there exist
infinitely many arithmetical progression of three different primes
$p$, such that $p+2=\mathrm{P}_2$. (A weaker result of this type
was previously obtained by Tolev \cite{Tol2}). In the paper 
\cite{GrT2} Green and Tao state that using the same method their result can be
extended for progression of $k$ terms, where $k$ is arbitrary
large.

In the present paper we consider another popular problem with
primes and study it with primes of the form specified above.

Let $\alpha $ be irrational real number, $\beta $ be real and let
$||x||=~\min\limits_{n\in \mathbb{Z}}|x-~n|$. In 1947 Vinogradov \cite{Vin4} proved that if
$0<\theta<1/5$ then there are infinitely many primes
$p$ such that
\begin{equation*}
    ||\alpha p+\beta||<p^{-\theta }.
\end{equation*}
Latter the upper bound for $\theta $ was improved and the
strongest published result is due to Heath-Brown and Jia
\cite{Heath2} with $\theta<16/49$. 
We shall prove the following

\begin{theorem} Let $\alpha \in \mathbb{R}\backslash \mathbb{Q}$ ,
$\beta \in \mathbb{R}$
 and let $0<~\theta\le ~1/100$. Then there are
infinitely many primes $p$ satisfying $p+2=\mathrm{P}_4$ and such
that
\begin{equation}\label{alfabeta}
    ||\alpha p+\beta||<p^{-\theta }.
\end{equation}
\end{theorem}

Other versions of this theorem are also possible, but our
intention is to present here a result with $r$ as small as possible and
for this $r$ to find some (not necessarily the biggest possible)
$\theta $.

\bigskip

{\bf Acknowledgements:} The first author was supported by the Sofia University Scientific Fund, grant 055/2007. The second author was supported by the Ministry of Science and Education of Bulgaria, grant MI-1504/2005.

\section{Notation}

Let $N$ be a sufficiently large real number and
$\delta,\,\eta,\,\theta,\,\rho,\,\kappa$ be real constants
satisfying
\begin{equation}\label{us1}
    0<\theta<\eta<\frac{\delta}{2}<\frac{1}{4},\qquad
    \eta<\rho<\delta, \qquad 0<\kappa, \qquad 0 < \theta \le \frac{1}{100} .
\end{equation}
We shall specify $\delta, \eta, \rho$ and $\kappa$  latter. We put
\begin{equation}\label{us2}
    \begin{split}
    z&=N^{\eta },\;\;\, y=N^{\rho },\quad D=N^{\delta },\\
    \Delta &=\Delta (N)=N^{-\theta },\quad \,H=\Delta ^{-1}\log ^2N.
    \end{split}
\end{equation}
By $p$ and $q$ we always denote primes. As usual $\Omega (n),\,\varphi (n),\,\mu (n),\,\Lambda
(n)$ denote respectively the numbers of prime factors of $n$
counted with the multiplicity, Euler's function, M\"{o}bius'
function and Mangoldt's function. We denote by $\tau _k(n)$ the
number of solutions of the equation $m_1m_2\ldots m_k=n$ in
natural numbers $m_1,\,\ldots,m_k$ and $\tau _2(n)=\tau (n)$. Let
$(m_1,\,\ldots \,,m_k)$ and $[m_1,\,\ldots \,,m_k]$ be the
greatest common divisor and the least common multiple of
$m_1,\,\ldots,m_k$ respectively. Instead of $m\equiv n\,\pmod {k}$
we write for simplicity $m\equiv n(k)$. As usual, $[y]$ denotes
the integer part of $y$, $||y||$ -- the distance from $y$ to the
nearest integer, $e(y)=e^{2\pi iy}$. For positive $A$ and $B$ we
write $A\asymp B$ instead of $A\ll B\ll A$. The letter
$\varepsilon$ denotes an arbitrary small positive number, not the
same in all appearances. For example this convention allows us to
write $x^{\varepsilon }\log x\ll x^{\varepsilon }$.

\section{Proof of the theorem}
We take a periodic with period 1 function such that
\begin{equation}\label{hi}
\begin{aligned}
    0<\chi (t) &<1 \quad \mbox { if }\quad  -\Delta< t< \Delta;\\
    \chi (t) &=0 \quad \mbox { if }\quad\quad\;  \Delta \le t\le 1-\Delta,
\end{aligned}
\end{equation}
and which has a Fourier series
\begin{equation}\label{hi1}
    \chi (t)= \Delta +\sum\limits_{|k|>0}g(k)e(kt),
\end{equation}
with coefficients satisfying

\begin{align}\label{hi2}
    g(0)&=\Delta,&\nonumber\\
    g(k)&\ll \Delta \;\;\;\;\;\mbox{ for all } k,\\
    \sum\limits_{|k|>H}|g(k)|&\ll
    N^{-1}.\nonumber
\end{align}

The existence of such a function is a consequence of a well known
lemma of Vinogradov (see \cite{Kar}, ch. 1, \S 2).

Consider the sum
\begin{equation}\label{gama}
    \Gamma =\Gamma (N)= \sum\limits_{N/2<p\le N \atop{(p+2,\,P(z))=1}}\chi (\alpha p+\beta
)\,T_p\,\log p,
\end{equation}
where
\begin{equation}\label{prod}
    P(z)=\prod\limits_{2<p\le z}p
\end{equation}
and
\begin{equation}\label{Tp}
    T_p=1-\kappa\sum\limits_{z<q\le y\atop{q|p+2}}\bigg(1-\frac{\log q}{\log y}\bigg).
\end{equation}
Obviously
\begin{equation}\label{gama1}
    \Gamma (N)\le \Gamma _1
\end{equation}
where $\Gamma _1$ is the sum of the terms of $\Gamma (N)$ for
which $T_p>0$. Denote by $\Gamma _2$ the sum of the term of
$\Gamma _1$ for which $\mu ^2(p+2)=0$. It is clear that
\begin{align}\label{gama3}
    0 \le \Gamma _2 
    \ll \sum\limits_{z\le q}\sum\limits_{n\le N\atop{n+2\equiv
    0(q^2)}}\log n
    &
    \ll \log N \sum\limits_{z\le q\le
    \sqrt{N+2}}\bigg(\frac{N}{q^2}+1\bigg) \\
    &
    \ll \frac{N^{1+\varepsilon}}{z}+N^{\frac{1}{2}+\varepsilon}\ll
    N^{1-\eta +\varepsilon}. \notag
\end{align}
We also remove from $\Gamma _1$ the term ( if such exist ) for
which $N-~2<~p\le ~N$ and the resulting error is $O(\log N)$.
Therefore
\begin{equation}\label{gama2}
    \Gamma \le\Gamma _3+O(N^{1-\eta +\varepsilon}),
\end{equation}
where
\begin{equation*}
    \Gamma _3=\sum\chi (\alpha p+\beta)\,T_p\,\log p
\end{equation*}
and where the summation is taken over the primes $p$, satisfying
\begin{equation}\label{Tp1}
    N/2<p\le N-2,
\end{equation}
\begin{equation}\label{Tp2}
    T_p>0,\quad\mu ^2(p+2)=1, \quad (p+2,P(z))=1.
\end{equation}
Assume that
\begin{equation}\label{gama5}
    \Gamma (N)\gg \frac{\Delta N}{\log N}.
\end{equation}
Then from~\eqref{us1},~\eqref{us2} and~\eqref{gama2} we get
$\Gamma _3>0$. Hence there exist a prime $p$
satisfying~\eqref{Tp1},~\eqref{Tp2} and such that
\begin{equation}\label{Tp3}
   \chi (\alpha p+\beta)>0.
\end{equation}
From~\eqref{us2},~\eqref{hi},~\eqref{Tp1} and~\eqref{Tp3} it
follows that this prime satisfies~\eqref{alfabeta}.

On the other hand, from the properties of the weights $T_p$ (see
\cite{Hal}, ch.~9 ) it follows that if $p$
satisfies~\eqref{Tp1},~\eqref{Tp2} then
\begin{equation*}
    \Omega (p+2)\le \frac{1}{\kappa}+\frac{1}{\rho}.
\end{equation*}
We see that to prove our theorem it is enough to determine the
constants $\delta,\,\eta,\,\theta,\,\rho,\,\kappa $ in such a way
that:
\begin{description}
    \item[\textit{I}] There exist a sequence $\{N_j\}_{j=1}^{\infty}$
    such that
\[
    \lim\limits_{j\to\infty}N_j=\infty, \qquad
    \Gamma(N_j)\gg \frac{\Delta (N_j)\,N_j}{\log
    N_j}\,,\;\;\;j=1,\,2,\,3,\ldots
\]
    \item[\textit{II}] We have
    \begin{equation}\label{gama8}
    \frac{1}{\kappa}+\frac{1}{\rho}<5
\end{equation}
\end{description}

Using~\eqref{gama} and~\eqref{Tp} we write $\Gamma $ as
\begin{equation}\label{gamaPG}
    \Gamma =\Phi - \kappa G,
\end{equation}
where
\begin{equation}\label{phi}
    \Phi =\sum\limits_{N/2<p\le N \atop{(p+2,\,P(z))=1}}\chi (\alpha p+\beta )\log p
\end{equation}
and
\begin{equation}\label{g}
    G=\sum\limits_{N/2<p\le N \atop{(p+2,\,P(z))=1}}\chi (\alpha p+\beta )\log p\sum\limits_{z<q\le y\atop{q|p+2}}\bigg (
1-\frac{\log q}{\log y}\bigg).
\end{equation}
We shall estimate $\Phi$ from below and $G$ from above.

Consider the sum $\Phi$. We apply a lower bound linear sieve. We
take the lower Rosser weights $\lambda ^-(d)$ of order $D$. For
the definition and their properties we refer the reader to
\cite{Iwa1}, \cite{Iwa2}. In particular we shall use that the
Rosser weights are real numbers such that
\begin{equation}\label{lam1}
    |\lambda ^-(d)|\le 1,\;\;\;\lambda ^-(d)=0 \;\;\;\mbox{if}\;\;\; d>D \;\;\; \mbox{or}\;\;\; \mu ^2(d)=0,
\end{equation}
\begin{equation}\label{lam2}
\sum\limits_{d|A}\lambda ^-(d)\le \begin{cases}1,\;\;\mbox{  if  }A=1,\\
 0,\;\; \mbox{  if }A\in \mathbb{N},\;A>1.
 \end{cases}
\end{equation}
We shall also use that if
\begin{equation}\label{s}
    s=\frac{\log D}{\log z}=\frac{\delta}{\eta} \qquad\mbox{  and  }\qquad 2<s< 4
\end{equation}
then
\begin{equation}\label{sumf}
    \sum\limits_{d|P(z)}\frac{\lambda ^-(d)}{\varphi (d)}\ge
    \Pi (z)\Bigg(\frac{2e^{\gamma}\log (s -1)}{s}+ O\bigg((\log
N)^{-1/3}\bigg)\Bigg),
\end{equation}
where
\begin{equation}\label{piz}
    \Pi (z)=\prod\limits_{2<p\le z}\bigg(1-\frac{1}{p-1}\bigg).
\end{equation}
From this place onwards we assume that
\begin{equation}\label{deltaovereta}
    2<\frac{\delta}{\eta}< 4,
\end{equation}
so the inequality~\eqref{sumf} holds.
Using~\eqref{phi},~\eqref{lam2} we get
\begin{align}\label{phi0}
    \Phi\ge \Phi _1&=\sum\limits_{N/2<p\le N }\chi (\alpha p+\beta )\log p
    \sum\limits_{d|(p+2,\,P(z))}\lambda ^{-}(d)\\
&=\sum\limits_{d|P(z)}\lambda ^{-}(d) \sum\limits_{N/2<p\le N
\atop{p+2\equiv 0\,(d)}}\chi (\alpha p+\beta )\log p.\nonumber
\end{align}
Form~\eqref{hi1},~\eqref{hi2} we find that
\begin{equation}\label{phi1}
    \Phi _1=\Delta (\Phi _2+\Phi_3)+O(1),
\end{equation}
where
\begin{eqnarray}
    \Phi _2&=&\sum\limits_{d|P(z)}\lambda ^-(d)
    \;
    \sum\limits_{N/2<p\le N
\atop{p+2\equiv 0\,(d)}}\log p,\nonumber\\
\label{phi3expres}
    \Phi _3&=&\sum\limits_{d|P(z)}\lambda ^-(d)
    \;
    \sum\limits_{0<|k|\le H}c(k)
    \;
    \sum\limits_{N/2<p\le N \atop{p+2\equiv 0\,(d)}}e(\alpha pk)\log
    p,\\
\label{ck}
    c(k)&=&\Delta ^{-1}g(k)e(\beta k)\ll 1.
\end{eqnarray}

Consider $\Phi _2$. From Bombieri-Vinogradov's theorem (see
\cite{Dav}, ch.~24), \eqref{us1},~\eqref{us2},~\eqref{lam1} we see
that
\begin{equation}\label{phi2}
    \Phi _2=\frac{N}{2}\sum\limits_{d|P(z)}\frac{\lambda ^{-}(d)}{\varphi (d)}+O\bigg(\frac{N}{(\log
    N)^2}\bigg).
\end{equation}
It is well known that the product defined
by~\eqref{piz} satisfies
\begin{equation}\label{pez}
    \Pi (z)\asymp \frac{1}{\log z}.
\end{equation}
Therefore from~\eqref{us1},~\eqref{us2},~\eqref{sumf},~\eqref{phi2}, \eqref{pez} we find
that
\begin{equation}\label{phi3}
    \Phi _2\ge e^{\gamma}N\,\Pi (z)\,\frac{\log (s -1)}{s}
+O\bigg(\frac{N}{(\log N)^{4/3}}\bigg),
\end{equation}
where $s$ is specified by~\eqref{s}. We shall study the sum $\Phi _3$ later.

Consider now the sum $G$, defined by~\eqref{g}. We write it in the
form
\begin{equation}\label{G}
    G=\sum\limits_{z<q<y}\bigg(1-\frac{\log q}{\log y}\bigg)
    \sum\limits_{N/2<p\le N\atop{p+2\equiv 0\, (q)\atop{(p+2,\,P(z))=1}}}\chi (\alpha p+\beta )\log p
\end{equation}
and then apply an upper bound linear sieve. Let $\lambda _q(d)$ be
the upper Rosser weights of order $\ds\frac{D}{q}$. We know that
\begin{equation}\label{lamq1}
    |\lambda _q(d)|\le 1,\;\;\;\lambda _q(d)=0\;\;\;\mbox{if}\;\;\; d>\frac{D}{q}\;\;\;\mbox{or}\;\;\;\mu ^2(d)=0,
\end{equation}
\begin{equation}\label{lamq2}
\sum\limits_{d|A}\lambda _q(d)\ge \begin{cases}1,\;\; \mbox{  if  }A=1,\\
 0,\;\; \mbox{  if }A\in \mathbb{N},\;A>1.
 \end{cases}
\end{equation}
We shall also use that if
\begin{equation}\label{dy}
    s_1=\frac{\log (D/q)}{\log z}\quad \mbox{and} \quad 1<s_1<3
\end{equation}
then
\begin{equation}\label{sumF}
    \sum\limits_{d|P(z)}\frac{\lambda _q(d)}{\varphi (d)}\le
\Pi(z)\Bigg(\frac{2e^{\gamma}}{s_1}+ O((\log N)^{-1/3})\Bigg)
\end{equation}
From this place onwards we assume that
\begin{equation}\label{etaro}
    \eta +\rho< \delta.
\end{equation}
Then using also~\eqref{deltaovereta} we see that the
condition~\eqref{dy} holds, consequently~\eqref{sumF} is true.
From~\eqref{G} --~\eqref{lamq2} we find
\begin{eqnarray}\label{gam0}
    &G\le &G_1=\\
    &\quad =&\sum\limits_{z<q<y}\bigg(1-\frac{\log q}{\log y}\bigg)\sum\limits_{N/2<p\le N\atop{p+2\equiv 0\,(q)}}
\chi (\alpha p+\beta )\log p\sum\limits_{d|(p+2,\,P(z))}\lambda _q(d)\nonumber\\
    &\quad =&\sum\limits_{z<q<y}\bigg(1-\frac{\log q}{\log y}\bigg)
\sum\limits_{d|P(z)}\lambda _q(d)\sum\limits_{N/2<p\le
N\atop{p+2\equiv 0\,(qd)}}
\chi (\alpha p+\beta )\log p\nonumber\\
&\quad =&\sum\limits_{m\le D}\gamma (m)\sum\limits_{N/2<p\le
N\atop{p+2\equiv 0\,(m)}} \chi (\alpha p+\beta )\log p,\nonumber
\end{eqnarray}
where
\begin{equation}\label{gamm}
    \gamma (m)=\sum\limits_{z<q<y\atop{d|P(z)\atop{qd=m}}}\bigg(1-\frac{\log q}{\log y}\bigg)\lambda _q(d).
\end{equation}
Using~\eqref{prod},~\eqref{lamq1} and~\eqref{gamm} we easily see
that
\begin{equation}\label{gamm1}
    |\gamma (m)|\le 1.
\end{equation}
From~\eqref{hi1},~\eqref{hi2} and~\eqref{gam0} we find
\begin{equation}\label{g1}
    G_1=\Delta (G_2+G_3)+O(1),
\end{equation}
where
\begin{align}
    G_2&=\sum\limits_{m\le D}\gamma (m)\sum\limits_{N/2<p\le N\atop{p+2\equiv 0\,(m)}}
\log p,\nonumber\\
\label{ge3} G_3&=\sum\limits_{m\le D}\gamma
(m)\sum\limits_{0<|k|\le H}c(k) \sum\limits_{N/2<p\le
N\atop{p+2\equiv 0\,(m)}} e(\alpha pk)\log p,
\end{align}
and where $c(k)$ satisfies~\eqref{ck}.

We apply again Bombieri-Vinogradov's theorem
and~\eqref{us1},~\eqref{us2}, \eqref{gamm1} to find that
\begin{equation}\label{g2}
    G_2=\frac{N}{2}\sum\limits_{m\le D}\frac{\gamma (m)}{\varphi (m)}+O\bigg(\frac{N}{(\log N)^2}\bigg).
\end{equation}
Using~\eqref{sumF},~\eqref{gamm} we obtain
\begin{align}\label{summ}
\sum\limits_{m\le D}\frac{\gamma (m)}{\varphi
(m)}&=\sum\limits_{z<q<y}\bigg(1-\frac{\log q}{\log y}\bigg)
\sum\limits_{d|P(z)}\frac{\lambda _{q}(d)}{\varphi (q d)}\\
&=\sum\limits_{z<q<y}\bigg(1-\frac{\log q}{\log
y}\bigg)\frac{1}{q-1}
\sum\limits_{d|P(z)}\frac{\lambda _{q}(d)}{\varphi (d)}\nonumber \\
&\le\sum\limits_{z<q<y}\bigg(1-\frac{\log q}{\log
y}\bigg)\frac{1}{q-1} \nonumber\\
&\quad\times\Pi (z)\Bigg(2e^{\gamma}\bigg(\frac{\log (D/q)}{\log
z}\bigg)^{-1}+ O\bigg((\log N)^{-1/3}\bigg)\Bigg).\nonumber
 \end{align}
Therefore from~\eqref{us2},~\eqref{pez},~\eqref{g2},~\eqref{summ}
we get
\begin{equation}\label{g2in}
    G_2\le e^{\gamma}N\,\Pi (z)
\sum\limits_{z<q<y}\bigg(1-\frac{\log q}{\log
y}\bigg)\frac{1}{q-1} \bigg(\frac{\log (D/q)}{\log
z}\bigg)^{-1}+O\bigg(\frac{N}{(\log N)^{4/3}}\bigg).
\end{equation}

Now we are in a position to find a lower bound for the sum $\Gamma$.
From~\eqref{gamaPG},~\eqref{phi0},~\eqref{phi1},~\eqref{phi3},~\eqref{gam0},~\eqref{g1},~\eqref{g2in}
it follows that
\begin{equation}\label{gam1}
    \Gamma\ge e^{\gamma}\Delta \,N\,\Pi (z)\,\Sigma +
O\bigg(\frac{\Delta N}{(\log N)^{4/3}}\bigg)+O\bigg(\Delta |\Phi
_3-\kappa G_3|\bigg),
\end{equation}
where
\begin{equation}\label{gam2}
    \Sigma =\frac{\log (s-1)}{s}-\kappa
\sum\limits_{z<q<y}\bigg(1-\frac{\log q}{\log
y}\bigg)\frac{1}{q-1}\bigg(\frac{\log (D/q)}{\log z}\bigg)^{-1}
\end{equation}
and where $s$ is specified by~\eqref{s}. Using partial summation
and the prime number theorem it is easy to prove that
\begin{equation}\label{gam3}
    \Sigma=\Sigma_0+O\bigg(\frac{1}{\log N}\bigg),
\end{equation}
where
\begin{equation}\label{sigma}
    \Sigma_0=\frac{\log (s -1)}{s}-\kappa \eta
\int\limits_{\eta}^{\rho}\bigg(\frac{1}{u}-\frac{1}{\rho}\bigg)
    \frac{1}{\delta-u}\,du.
\end{equation}
Therefore, using~\eqref{us2},~\eqref{pez},~\eqref{gam1} we get
\begin{equation*}
   \Gamma\ge e^{\gamma}N\,\Delta \,\Pi (z)\,\Sigma _0 +
O\bigg(\frac{\Delta N}{(\log N)^{4/3}}\bigg)+O\bigg(\Delta |\Phi
_3-\kappa G_3|\bigg).
\end{equation*}

Now we shall see that if $N$ is a term of a suitable sequence
tending to infinity then the last error term in the formula above can be omitted.
The following lemma holds:

\begin{lemma} Suppose $\alpha \in \mathbb{R}\backslash \mathbb{Q}$ and
\begin{equation}\label{deltaeta}
   \delta +\theta <\frac{1}{3}\,.
\end{equation}
Let $\xi(d),\,c(k)$ be complex numbers defined for $d\le
D,\;0<|k|\le H$, where $D$ and $H$ are specified by~\eqref{us2}, and let
\begin{equation} \label{25oct1}
    \xi(d)\ll 1,\qquad c(k)\ll 1.
\end{equation}
Then there exist
a sequence $\{N_j\}_{j=1}^{\infty},\;\lim\limits_{j\to
\infty}N_j=\infty$, such that if
\begin{equation}\label{sn1}
   S(N)= \sum\limits_{d\le D}\xi(d)
    \sum\limits_{1\le |k|\le H}c(k)
    \sum\limits_{N/2<p\le N\atop{p+2\equiv 0\,(d)}}e(\alpha pk)\log p
\end{equation}
then we have
\begin{equation*}
    S(N_j)\ll \frac{N_j}{\log ^2N_j}\,,\;\;\;\;j=1,\,2,\,3,\ldots\;.
\end{equation*}
\end{lemma}

We will present the proof of this Lemma in the next section.

\bigskip

From~\eqref{phi3expres},~\eqref{ge3} we see that the quantity
$\Phi_3-\kappa G_3$ can be written as a sum of type~\eqref{sn1}
with $\xi (d)=\lambda ^*(d)-\kappa\gamma (d)$, where 
$\lambda ^* (d) =\lambda ^-(d)$ if $d|P(z)$ and 
$\lambda ^* (d) =0$ otherwise. Using our Lemma we
see that there exist a sequence tending to infinity such that if
$N$ is its term then
\begin{equation}\label{gam7}
   \Gamma\ge e^{\gamma}\Delta \,N\,\Pi (z)\,\Sigma _0 +
O\bigg(\frac{\Delta N}{(\log N)^{4/3}}\bigg).
\end{equation}

We put
\begin{equation*}
    \rho=0.23,\quad\delta=0.315,\quad\eta=0.08,\quad\kappa=1.58\,.
\end{equation*}
Then it is easy to verify that the
conditions~\eqref{us1},~\eqref{gama8},~\eqref{deltaovereta},~\eqref{etaro},
~\eqref{deltaeta} are fulfilled  and also 
\[
\Sigma_0>0 .
\]
From
the last inequality and~\eqref{us2},~\eqref{pez},~\eqref{gam7} it
follows that there exist a constant $c>0$ such that
for any $N$ from our sequence we have
\begin{equation*}
        \Gamma\ge c\,\frac{\Delta N}{\log N}>0.
\end{equation*}
This completes the proof of the theorem.

\section{Proof of Lemma}

Since $\alpha$ is irrational we see, using Dirichlet's theorem,  that there are infinitely many integers $A$ and natural numbers $Q$ such that
\begin{equation*}
    \bigg|\alpha - \frac{A}{Q}\bigg|<\frac{1}{Q^2} .
\end{equation*}
For any such $Q$ we choose $N$ in a suitable way
(see \eqref{kraj} ) and in this way define our sequence $\{N_j\}_{j=1}^{\infty}$.

It is clear that
\begin{equation}\label{sn}
    S(N)=W+O(HN^{\frac{1}{2}+\varepsilon }),
\end{equation}
where
\begin{equation*}
    W=\sum\limits_{N/2<n\le N}\Lambda (n)
    \sum\limits_{1\le |k|\le H}c(k)e(\alpha nk)
    \sum\limits_{d\le D\atop{d|n+2\atop{2\nmid d}}}\xi(d).
\end{equation*}
Using Heath-Brown's identity \cite{Heath} with parameters
\begin{equation}\label{uvw}
P=N/2,\,P_1=N,\,u=2^{-7}N^{\frac{\delta}{2}},\,v=2^7N^{\frac{1}{3}},\,w=N^{\frac{1}{2}-\frac{\delta}{4}}.
\end{equation}
we decompose the sum $W$ as a linear combination of $O(\log ^6N)$
sums of first and second type. The sums of the first type are
\begin{equation*}
    W_1=\sum\limits_{M<m\le M_1}a_m
    \sum\limits_{L<l\le L_1\atop{N/2<ml\le N}}
    \sum\limits_{0<|k|\le H}c(k)e(\alpha mlk)
    \sum\limits_{d\le D\atop{d|ml+2\atop{2\nmid d}}}\xi(d)
\end{equation*}
and
\begin{equation*}
    W'_1=\sum\limits_{M<m\le M_1}a_m
    \sum\limits_{L<l\le L_1\atop{N/2<ml\le N}}\log l
    \sum\limits_{0<|k|\le H}c(k)e(\alpha mlk)
    \sum\limits_{d\le D\atop{d|ml+2\atop{2\nmid d}}}\xi(d) ,
\end{equation*}
where
\begin{equation}\label{Usl1}
M_1\le 2M, \quad L_1\le 2L, \quad   ML\asymp N, \quad  L\ge w, \quad  a_m \ll  N^{\varepsilon}.
\end{equation}

The sums of the second type are
\begin{equation*}
    W_2=\sum\limits_{M<m\le M_1}a_m
    \sum\limits_{L<l\le L_1\atop{N/2<ml\le N}}b_l
    \sum\limits_{0<|k|\le H}c(k)e(\alpha mlk)
    \sum\limits_{d\le D\atop{d|ml+2\atop{2\nmid d}}}\xi(d) ,
\end{equation*}
where
\begin{equation}\label{Usl3}
M_1\le 2M, \quad L_1\le 2L, \quad ML\asymp N, \quad u\le L\le v, \quad a_m, \; b_l \ll  N^{\varepsilon}.
\end{equation}

First we estimate the sums of second type. We have
\begin{equation*}
    W_2 \ll N^{\varepsilon }\sum\limits_{M<m\le M_1}
   \bigg | \sum\limits_{L<l\le L_1\atop{N/2<ml\le N}}b_l
    \sum\limits_{0<|k|\le H}c(k)e(\alpha mlk)
    \sum\limits_{d\le D\atop{d|ml+2\atop{2\nmid d}}}\xi(d)\bigg |.
\end{equation*}
Applying the Cauchy inequality and \eqref{25oct1}, \eqref{Usl3} we get
\begin{align*}
    |W_2|^2&\ll N^{\varepsilon }M\sum\limits_{M<m\le M_1}
   \bigg | \sum\limits_{L<l\le L_1\atop{N/2<ml\le N_1}}b_l
    \sum\limits_{0<|k|\le H}c(k)e(\alpha mlk)
    \sum\limits_{d\le D\atop{d|ml+2\atop{2\nmid d}}}\xi (d)\bigg |^2\\
    &\ll N^{\varepsilon }M
    \sum\limits_{\substack{d_1,\,d_2\le D \\ 2 \nmid{d_1d_2} }} \;
    \sum\limits_{0<k_1,\,k_2\le H} \;
    \sum\limits_{L<l_1,\,l_2\le L_1}|V|,
\end{align*}
where
\begin{gather*}
   V= \sum\limits_{M'<m\le M'_1\atop{l_im+2\equiv 0(d_i),\,i=1,2}}e(\alpha
    m(k_1l_1-k_2l_2)),\\
     M'=\max \bigg \{ \frac{N}{2l_1},\,\frac{N}{2l_2},\,M\bigg \},\;
    M'_1=\min \bigg \{ \frac{N}{l_1},\,\frac{N}{l_2},\,M_1\bigg \}.
\end{gather*}
If the system of congruences
\begin{equation}\label{cond}
   \left |
\begin{array}
    ll_1m+2\equiv 0(d_1)\\
    l_2m+2\equiv 0(d_2).
\end{array}
\right.
\end{equation}
has no solution then $V=0$. Assume that the system \eqref{cond}
has a solution. Then there exist an integer
$f=f(l_1,\,l_2,\,d_1,\,d_2)$ such that \eqref{cond} is equivalent to $m\equiv f([d_1,\,d_2])$
and therefore
\begin{align*}
    V&=\sum\limits_{M'<m\le M'_1\atop{m\equiv f([d_1,\,d_2])}}e(\alpha
    m(k_1l_1-k_2l_2))\\
    &=e(\alpha f(k_1l_1-k_2l_2))
    \sum\limits_{\frac{M'-f}{[d_1,\,d_2]}<s\le \frac{M'_1-f}{[d_1,\,d_2]}}e(\alpha s[d_1,\,d_2](k_1l_1-k_2l_2)).
\end{align*}
From~\eqref{us2},~\eqref{deltaeta},~\eqref{uvw},~\eqref{Usl3} it
follows that
\begin{equation}\label{MD}
    M\gg \frac{N}{v}\gg D^2 .
\end{equation}
Applying Lemma 4 from \cite{Kar},
ch. 6, \S 2, we get
\begin{equation}\label{v}
    V\ll \begin{cases}\ds \frac{M}{[d_1,\,d_2]}\, ,&\mbox{  if }k_1l_1=k_2l_2,\\[10pt]
 \ds \min \bigg \{\frac{M}{[d_1,\,d_2]},\,\frac{1}{||\alpha
 (k_1l_1-k_2l_2)[d_1,\,d_2]||}\bigg \}\,
 ,&\mbox{  if }k_1l_1\neq k_2l_2.
\end{cases}
\end{equation}

Therefore
\begin{equation}\label{w21}
    |W_2|^2\ll N^{\varepsilon}M\bigg(MV_1+V_2\bigg),
\end{equation}
where
\begin{align*}
    V_1 &=\sum\limits_{d_1,\,d_2\le D}\frac{1}{[d_1,\,d_2]}
      \;
    \sum\limits_{0<k_1,\,k_2\le H}
      \;
    \sum\limits_{L<l_1,\,l_2\le L_1\atop{k_1l_1=k_2l_2}}1,\\
    V_2 &=\sum\limits_{d_1,\,d_2\le D}
      \;
    \sum\limits_{0<k_1,\,k_2\le H}
      \;
    \sum\limits_{L<l_1,\,l_2\le
    L_1\atop{k_1l_1\ne k_2l_2}}
    \min \bigg \{\frac{M}{[d_1,\,d_2]},\,\frac{1}{||\alpha
 (k_1l_1-k_2l_2)[d_1,\,d_2]||}\bigg \}.
\end{align*}

It is clear that
\begin{multline}\label{v1}
    V_1 \ll \sum\limits_{h\le
    D^2}\frac{1}{h}\sum\limits_{[d_1,\,d_2]=h}1
    \sum\limits_{n\le 2HL}\tau ^2(n)\ll
    N^{\varepsilon}HL
    \sum\limits_{h\le D^2}\frac{\tau ^2(h)}{h}\ll
    N^{\varepsilon}HL.
\end{multline}
Consider $V_2$ we have
\begin{align*}
    V_2 &\ll \sum\limits_{h\le D^2}\tau ^2(h)
    \sum\limits_{0<|r|\le 2HL}
    \min \bigg \{\frac{M}{h},\,\frac{1}{||\alpha rh||}\bigg \}
    \sum\limits_{0<n_1,\,n_2\le 2HL\atop{n_1-n_2=r}}\tau (n_1)\tau (n_2)\\
    &\ll N^{\varepsilon}HL\sum\limits_{h\le D^2}
      \;
    \sum\limits_{0< r \le 2HL}
    \min \bigg \{\frac{M}{h},\,\frac{1}{||\alpha rh||}\bigg \}\\
    &\ll N^{\varepsilon}HL \sum\limits_{m\le 2D^2HL}
    \min \bigg \{\frac{HLM}{m},\,\frac{1}{||\alpha m||}\bigg \}.
\end{align*}
Since $M \gg D^2$ (see~\eqref{MD}) we can apply Lemma 2.2 from
\cite{Von}, ch. 2, \S 2.1 and get
\begin{equation}\label{v2}
   V_2\ll N^{\varepsilon }\bigg(\frac{H^2L^2M}{Q}+D^2H^2L^2+HLQ\bigg).
\end{equation}
From~\eqref{Usl3},~\eqref{w21} -- \eqref{v2} we
obtain
\begin{align*}
    |W_2|^2&\ll N^{\varepsilon }
    \bigg(HLM^2+\frac{H^2L^2M^2}{Q}+D^2H^2L^2M+HLMQ\bigg)\\
    &\ll N^{\varepsilon }\bigg(\frac{HN^{2}}{u}+\frac{H^2N^{2}}{Q}+
    D^2H^2Nv+HNQ\bigg).
\end{align*}
Hence
\begin{equation}\label{w2}
    W_2 \ll N^{\varepsilon }\bigg(\frac{H^{\frac{1}{2}}N}{u^{\frac{1}{2}}}+\frac{HN}{Q^{\frac{1}{2}}}+
    DHN^{\frac{1}{2}}v^{\frac{1}{2}}+H^{\frac{1}{2}}N^{\frac{1}{2}}Q^{\frac{1}{2}}\bigg).
\end{equation}

Now we shall estimate the sums of the first type. Using \eqref{25oct1}, \eqref{Usl1} we
get
\begin{equation}\label{w11}
    W_1 \ll N^{\varepsilon }
    \sum\limits_{d\le D\atop{2\nmid d}}
      \;
    \sum\limits_{0<k\le H}
      \;
    \sum\limits_{M<m\le M_1}|U|,
\end{equation}
where
\begin{gather}\label{u}
    U=\sum\limits_{L'<l\le L_1'\atop{ml+2\equiv 0(d)}}e(\alpha kml),\\
    L'=\max\bigg \{L,\,\frac{N}{2m}\bigg \},\;\;L'_1=\min\bigg \{L_1,\,\frac{N}{m}\bigg
    \}.\nonumber
\end{gather}
If $(m,\,d)>1$ then the sum $U$ is empty. Suppose now that
$(m,\,d)=1$. 
Then the congruence $ml+2 \equiv 0 (d)$
is equivalent to $l \equiv l_0 (d)$ for some integer
$l_0 = l_0(m, d)$. Hence we may write $U$ in the form
\begin{equation*}
    U=e(\alpha kml_0)
    \sum\limits_{\frac{L'-l_0}{d}< s\le \frac{L_1'+l_0}{d}}
    e(\alpha k m s d).
\end{equation*}
Using Lemma 4 from
\cite{Kar}, ch. 6, \S 2 we get
\begin{equation*}
    U\ll\min \bigg \{\frac{N}{md},\,\frac{1}{||\alpha kmd||}\bigg\} ,
\end{equation*}
consequently
\begin{align*}
    W_1&\ll N^{\varepsilon }
    \sum\limits_{d\le D}
      \;
    \sum\limits_{k\le H}
      \;
    \sum\limits_{M<m\le M_1}
    \min \bigg \{\frac{N}{md},\,\frac{1}{||\alpha kmd||}\bigg \}\\
    &\ll N^{\varepsilon }
    \sum\limits_{n\le 2MD}
      \;
    \sum\limits_{k\le H}
    \min \bigg \{\frac{N}{n},\,\frac{1}{||\alpha kn||}\bigg \}\\
    &\ll N^{\varepsilon}\sum\limits_{n\le 2MD}
     \;   
    \sum\limits_{k\le H}
    \min \bigg \{\frac{HN}{kn},\,\frac{1}{||\alpha kn||}\bigg \}\\
    &\ll N^{\varepsilon}
    \sum\limits_{s\le 2MDH}
    \min \bigg \{\frac{NH}{s},\,\frac{1}{||\alpha s||}\bigg \}.
\end{align*}
Using~\eqref{deltaeta},~\eqref{uvw},~\eqref{Usl1} we see that we
may apply again Lemma 2.2, \cite{Von}, ch. 2, \S 2.1 and we get
\begin{equation}\label{w1}
    W_1\ll N^{\varepsilon}\bigg(\frac{HN}{Q}+\frac{DHN}{w}+Q\bigg).
\end{equation}
We consider the sum $W_1'$ in the same manner and we find
\begin{equation}\label{w1'}
    W_1'\ll N^{\varepsilon}\bigg(\frac{HN}{Q}+\frac{DHN}{w}+Q\bigg).
\end{equation}
From~\eqref{w2},~\eqref{w1} and~\eqref{w1'} we obtain
\begin{multline*}
    W\ll N^{\varepsilon}\bigg(\frac{H^{\frac{1}{2}}N}{u^{\frac{1}{2}}}+\frac{HN}{Q^{\frac{1}{2}}}+
    DHN^{\frac{1}{2}}v^{\frac{1}{2}}+H^{\frac{1}{2}}N^{\frac{1}{2}}Q^{\frac{1}{2}}+
    \frac{DHN}{w}+Q\bigg).
\end{multline*}
We choose 
\begin{equation}\label{kraj}
N=Q^{\frac{2}{1+\theta}}
\end{equation}
and having in mind
\eqref{us1}, \eqref{deltaeta}, \eqref{uvw} we obtain
\begin{equation*}
    W\ll
    N^{1+\frac{\theta}{2}-\frac{\delta}{4}}+ N^{\frac{3(1+\theta)}{4}}+
    N^{\frac{2}{3}+\delta+\theta}
     \ll N^{1- \varpi} 
\end{equation*}
for some small constant $\varpi > 0 $. This proves our Lemma.

\bigskip
\bigskip

\vbox{
\hbox{Faculty of Mathematics and Informatics}
\hbox{Sofia University ``St. Kl. Ohridsky''}
\hbox{5 J.Bourchier, 1164 Sofia, Bulgaria}
\hbox{ }
\hbox{tlt@fmi.uni-sofia.bg}
\hbox{dtolev@fmi.uni-sofia.bg}}

\end{document}